\documentclass[12pt]{amsart}
\usepackage{amsmath}
\usepackage{amssymb,amscd}
\usepackage{amsfonts}
\usepackage{euscript}

\usepackage{graphicx}
\usepackage{color}

\DeclareMathAlphabet{\mathpzc}{OT1}{pzc}{m}{it}

\numberwithin{equation}{section}

\newtheorem{theorem}{Theorem}[section]

\newtheorem{corollary}[theorem]{Corollary}

\newtheorem{lemma}[theorem]{Lemma}
\newtheorem{prop}[theorem]{Proposition}

\theoremstyle{definition}
\newtheorem{defn}[theorem]{Definition}

\newtheorem{example}[theorem]{Example}
\newtheorem{ques}{Question}[section]
\newtheorem{remark}[theorem]{Remark}

\def \begineq{\begin{equation}}
\def \endeq{\end{equation}}

\def \bb{\mathbb}

\def \RR{{\bb{R}}}

\def \ZZ{{\bb{Z}}}

\def \({\left(}
\def \){\right)}
\def \<{\langle}
\def \>{\rangle}
\def \bar{\overline}

\begin{document}

\title[Orientability and equivariant cobordism of 2-torus manifolds
]{Orientability and equivariant oriented cobordism of 2-torus manifolds}

\author[S. Sarkar]{Soumen Sarkar}

\address{Department of Mathematics and Statistics, University of Calgary, Calgary, Canada T2N 1N4.}

\email{soumensarkar20@gmail.com}

%
%



\abstract
We give a necessary and sufficient condition for the orientability of a locally standard 2-torus manifold with a fixed point which generalizes previous results of Nakayama-Nishimura in 2005 and Soprunova-Sottile in 2013. We construct manifolds with boundary where the boundary is a disjoint union of locally standard 2-torus manifolds. We discuss equivariant oriented cobordism class of locally standard 2-torus manifolds.
\endabstract

\maketitle

\section{Introduction}\label{intro}
Group actions on topological spaces have many rich structures. Nice exposition of this subject can be found in Bredon's book \cite{Bre}. The standard faithful representation of $\ZZ_2^n$ on $\RR^n$ is determined by the reflections along coordinate hyper planes. Observing this phenomena Davis-Januskiewicz introduced the category of small covers in their pioneering paper \cite{DJ}. Generalizing the concept of small covers, L\"{u} and Masuda introduced the notion of 2-torus manifolds in \cite{LM}. Briefly, A $2$-torus manifold is a closed smooth manifold of dimension $n$ with an
effective action of a $2$-torus group $\ZZ_2^n$ of rank $n$. A 2-torus manifold is called locally standard if it is locally isomorphic to a faithful representation of $\ZZ_2^n$ on $\RR^n$. The orbit space $Q$ of an $n$-dimensional locally standard $2$-torus manifold $N$ by the $\ZZ_2^n$action is an $n$-dimensional nice manifold with corners. When $Q$ is a simple convex polytope, the manifold $N$ is called a small cover and studied its various topological properties in \cite{DJ}. A class of
examples of a small covers are real projective spaces $\mathbb{R} P^n$ and its products with standard actions of $\ZZ_2^n$. In this article, we study orientability and equivariant oriented cobordism of 2-torus manifolds. Our orientability results of a locally standard 2-torus manifold with a fixed point  generalizes the previous results of Nakayama-Nishimura \cite{NN}.

Cobordism was introduced by Lev Pontryagin in geometric work on manifolds, \cite{Pon}. In the early 1950's Ren\'{e} Thom \cite{Tho} showed that cobordism groups for the category of manifolds could be computed by results of homotopy theory using the Thom complex construction. The results of Thom can not be generalized to compute the equivariant cobordism groups. One of the reason is that Thom transversality may not be valid for equivariant category. There are several developments, such as \cite{Stong, Stong2, Mad, LT}, to study the equivariant cobordism groups for finite groups. In Section 1 of \cite{Qui}, Quillen discussed geometric interpretation of complex cobordism rings. Following his definition, we define the equivariant oriented cobordism groups and rings for the category of $2$-torus manifolds. 

Let $OS_n$ be the set of equivariant oriented cobordism classes of all $n$-dimensional 2-torus manifolds. Then $OS_n$ forms an abelian group with respect to disjoint union. Inspired by the work of Lu and Tan  such as \cite{Lu, LT}, we may ask the following question.
\begin{ques}\label{ques}
Does a class of $OS_n$ contain a small cover as its representative?
\end{ques}

The paper is organized as follows. In Subsections \ref{axdif} and \ref{condif}, we discuss the axiomatic and constructive definitions of locally standard 2-torus manifolds following \cite{LM}.
 In Subsection \ref{ori_2torus}, we give a necessary and sufficient condition for orientability of of 2-torus manifolds generalizing the Theorem 1.7 of Nakayama-Nishimura \cite{NN} and the Theorem 3.1 (for the smooth points of real toric variety) of Soprunova-Sottile \cite{SS2}. In Subsection \ref{eqi_con_sum} we briefly discuss equivariant connected sum of 2-torus manifolds. In Section \ref{def}, we construct manifolds with boundary where the boundary is a disjoint union of locally standard 2-torus manifolds. In Section \ref{cobor}, we discuss equivariant oriented cobordism class of locally standard 2-torus manifolds and give some sufficient conditions to answer Question \ref{ques}.

\section{Locally  standard 2-torus manifolds}\label{2torm}
Definition and equivariant classification of locally standard 2-torus manifolds are extensively discussed in \cite{LM}.
Following them, we discuss the axiomatic and constructive definitions of locally standard 2-torus manifolds. We begin by recalling the definition of manifold with corners following Section 6 of \cite{Da}. Various properties of manifold with corners and maps between manifolds with corners are studied in \cite{Jo}

\begin{defn}
(a) A Hausdorff topological space $Q \subset \RR^m$ is called an $n$-dimensional manifold with corners if any point $q \in Q$ has a neighborhood homeomorphic to an open subset in the positive cone $\RR^n_{\geq 0}= \{(x_1, \ldots, x_n)\in \RR^n: x_1 \geq 0, \ldots, x_n \geq 0\}$, where $n \leq m$, preserving the codimension function. Codimension one faces of $Q$ are called facets, and the faces of dimension 0 are called the vertices. 

(b) A manifold with corners is called nice if a codimension-2 face is the connected component of the  intersection of a unique collection of two facets.
\end{defn}

An $n$-dimensional simple polytope is a convex polytope each of whose vertex is the intersection of exactly $n$ facets. So simple polytopes are nice manifold with corners.  Note that a codimension-$k$ face of a nice manifold with corners is a connected component of the intersection of a unique collection of $k$ many facets. Throughout this article all manifold with corners are nice and subset of some $\RR^m$.

\subsection{Axiomatic definition}\label{axdif}
In this subsection, we recall axiomatic definition of $2$-torus manifolds. We denote the quotient additive group $\ZZ/2\ZZ$ by $\ZZ_2$ throughout this paper. Let $\rho_s : \ZZ_2^n \times \RR^n \to \RR^n$ be the standard action. That is, the action is generated by the reflections along the coordinate hyper planes. So the orbit space of this action is the positive cone $\RR^n_{\geq 0}$.
\begin{defn}
A smooth, closed, $n$-dimensional manifold $N$ is called locally standard $2$-torus manifold if the followings hold.
\begin{enumerate}
\item There is an effective action $\rho :\ZZ_2^n \times N \to N$.
\item Every point $y \in N $ has a $\ZZ_2^n$-invariant open neighborhood $U_y \subset N$, that is $\rho (\ZZ_2^n \times U_y) = U_y$.
\item There exists a homeomorphism $\psi : U_y \to V$, where $V$ is a $\ZZ_2^n$-invariant (that is $\rho_s (\ZZ_2^n \times V) = V$) open subset of $\RR^n$. 
\item There exists an isomorphism $\delta_y : \ZZ_2^n \to \ZZ_2^n$ such that
$\psi( \rho (t, x)) = \rho_s(\delta_y (t), \psi(x))$ for all $(t,x) \in \ZZ_2^n \times U_y$.
\end{enumerate}
\end{defn}
Let $Q$ be the orbit space of the $\ZZ_2^n$-action on $N$ and $\pi : N \to Q$ be the orbit map. By the
above definition, the orbit space $Q$ is an $n$-dimensional manifold with corners. Locally standardness of the action ensures that fixed point set correspond bijectively to the number of vertices of $Q$. Small covers correspond to the case when the orbit space is homeomorphic, as manifold with corners, to an $n$-dimensional simple polytope. Note that the map $\psi : U_y \to V$ satisfying the conditions (3) and (4) in the above definition is called a $\ZZ_2^n$-equivariant map. In the definition of locally standard 2-torus manifolds we did not assume that they are connected like in \cite{LM}.  

\begin{defn}
 A connected, closed, codimension-one submanifold of $N$ is called  $\ZZ_2$-$characteristic$ if it is a connected
component of the set which is fixed pointwise by a subgroup of $\ZZ_2^n$ isomorphic to $\ZZ_2$.
\end{defn}

Since $N$ is compact, $N$ has finitely many characteristic submanifolds. Note that the action of $\ZZ_2^n$ is free outside the union of all characteristic submanifolds. A connected component of the intersection of $k$ many $\ZZ_2$-characteristic submanifolds of N corresponds to a codimension-k face of $Q$. So a codimension-k face of $Q$ is a connected component of the intersection of $k$ many facets. That means, $Q$ is a nice manifold with corners.

\begin{example}\label{rpn}
The natural action of $\ZZ_2^n$ defined on the real projective space $\RR P^n$ by
\begin{equation}
(g_1, \ldots, g_n)\cdot [x_0, x_1, \ldots, x_n] \to [x_0, (-1)^{g_1}x_1, \ldots, (-1)^{g_n}x_n]
\end{equation}
is locally standard and the orbit space is homeomorphic as manifold with corners to the
standard $n$-simplex. Hence $\RR P^n$ is a small cover over the $n$-simplex $\Delta^n$.
\end{example}

\begin{example}
 Consider the sphere $S^{n} = \{(x_1, \ldots, x_n, x_{n+1}) \in \RR^{n+1} : x_1^2 + \cdots + x_{n+1} =1 \}$.
Let the action of $\ZZ_2^n$ on $S^n$ is generated by the reflection on $x_i=0$ planes for $i=1, \ldots, n$.
Then the action is locally standard and the orbit space is $\{(x_1, \ldots, x_n, x_{n+1}) \in S^n : x_1, \ldots,
x_n \geq 0\}$ which is not a simple polytope upto homeomorphism as manifold with corners. So $S^n$ is not a small cover but $2$-torus manifold.
\end{example}

\subsection{Definition by construction}\label{condif}
In this subsection, we discuss the definition of $2$-torus manifolds by construction.
Let $Q$ be an $n$-dimensional manifold with corners and $\mathcal{F}(Q)= \{F_i : i = 1, \ldots, m \}$ be the facets of $Q$.
\begin{defn}\label{z2char}
A function $\xi : \mathcal{F}(Q) \to \ZZ_2^n$ is called a $\ZZ_2$-characteristic function on $Q$ if 
the vectors $\{\xi(F_{i_1}), \ldots, \xi(F_{i_{\ell}})\}$ are linearly independent in $\ZZ_2^n$ whenever the intersection $F_{i_1}\cap \cdots \cap F_{i_{\ell}} $ is nonempty.
\end{defn}
The vector $\xi_i = \xi(F_i)$ is called the $\ZZ_2$-$dicharacteristic~ vector$ corresponding to the $i$-th facet $F_i$ of $Q$. Let $\tau : E \to Q$ be a principal $\ZZ_2^n$-bundle over $Q$. These two data will be referred to as {\it $\ZZ_2$-combinatorial model} and abbreviated as $((Q, \xi), (E, Q \tau))$. In Definition \ref{z2char}, we use `$\ZZ_2$-characteristic function' in contrast to `characteristic function' of \cite{LM} to avoid name conflict with the similar definition in other areas of toric topology.

\begin{remark}
\begin{enumerate}
\item If $n=1$ then $Q$ is either closed interval or a circle and if $n=2$ then $Q$ can be obtained from a closed surface by removing the interior points of a finite number of non-intersecting polygons, or polygons and eye shapes, or polygons and discs, or polygons and eye shapes and discs. For an example see Figure \ref{egsc3} (a). So, when $n \leq 2$, one can show that any $Q$ admits a $\ZZ_2$-characteristic function.

\item When $n = 3$, Four Color Theorem ensures that $Q$ admits a $\ZZ_2$-characteristic function if the boundary of $Q$ is a union of 2-spheres. 

\item  When $n \geq 4$, there is a simple convex polytope on which there is no $\ZZ_2$-characteristic function, see Nonexamples 1.22 in \cite{DJ}.
\end{enumerate}
\end{remark}

\begin{example}
The manifold with corners in Figure \ref{egsc3} $(a)$ is obtained by deleting the interior of the circle
$C$ and the triangle $v_4v_5v_6$ from the interior of the rectangle $v_0v_1v_2v_3$. The manifold with
corners in Figure \ref{egsc3} $(b)$ is obtained by deleting the interior of pentagon $v_0v_1v_2v_3v_4$
from the interior of the circle $C$. Some $\ZZ_2$-characteristic functions of these $2$-dimensional manifold
with corners is defined in the corresponding figure.
\begin{figure}[ht]
        \centerline{
           \scalebox{.76}{
            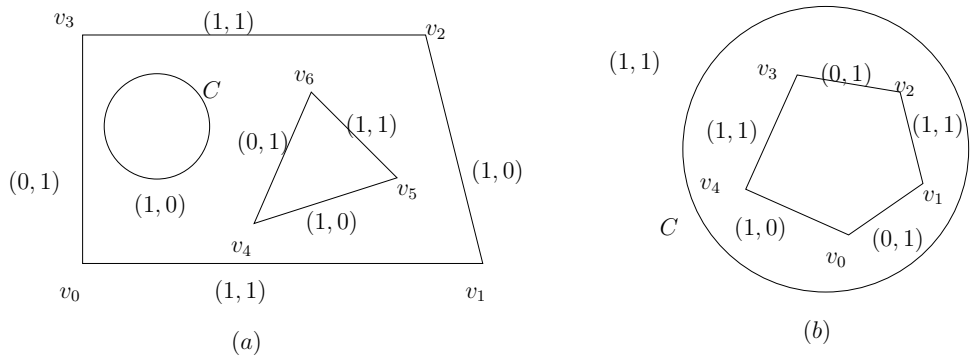
            }
          }
       \caption {Some $\ZZ_2$-characteristic functions.}
        \label{egsc3}
      \end{figure}
\end{example}

We also remark that a $\ZZ_2$-characteristic function arises naturally from a locally standard 2-torus manifold $N$. Let $N$ be a 2-torus manifold of dimension $n$ over the manifold $Q$ with corners. Note that a facet of $Q$ is the image of a $\ZZ_2$-characteristic function submanifold of $N$ by the quotient map $\pi : N \to Q$. To each element $F \in \mathcal{F}(Q)$ we assign the nonzero element of $\ZZ_2^n $ which fixes pointwise the characteristic submanifold $\pi^{-1}(F)$. The local standardness of $N$ implies that this assignment satisfies the linearly independent condition which is required to be a $\ZZ_2$-characteristic function $\xi$ on $Q$. For example, let $$N_i = \{[x_0, \ldots, x_n] \in \RR P^n : x_i =0\}.$$ Then the isotropy group of $N_i$ is the $i$-th factor subgroup of $\ZZ_2^n$ for $i= 1, \ldots, n$ and the isotropy group of $N_0$ is the diagonal subgroup of $\ZZ_2^n$. The corresponding $\ZZ_2$-characteristic function on $\Delta^n$ is called standard, since it is unique.

On the other hand, since $Q$ is a nice manifold with corners, the boundary $\partial{Q}$ has a color neighborhood in $Q$. That is, the complement $Q_c$ of this neighborhood is homeomorphic to $Q$ preserving the face structure. Note that $N_{\xi} = \pi^{-1}(Q_c)$ is the total space of a principal $\ZZ_2^n$-bundle, denoted by $\tau : E_N \to Q_c$. Since $Q_c$ is homeomorphic as manifold with corners to $Q$, we may assume $\tau : E_N \to Q$ is the associated principal $\ZZ_2^n$-bundle. Thus we
get a $\ZZ_2$-combinatorial data $((Q, \xi), (E_N, Q, \tau))$ from a locally standard 2-torus manifold $N$. 

Using the $\ZZ_2$-characteristic model $((Q, \xi), (E_N, Q, \tau)$, we can construct a locally standard 2-torus manifold $N(Q, \xi, \tau)$ similarly to the construction in Section 1.5 of \cite{DJ}. Since $Q$ is a manifold with corners, the space $E_N$ is a manifold with corners and each point $x \in Q$ belongs to the relative interior of a codimension-$\ell$ face $F$ which is a connected component of $F_{i_1} \cap \ldots \cap F_{i_{\ell}}$ for a unique collection of facets $F_{i_1}, \ldots, F_{i_{\ell}}$. Let $\ZZ_2(F)$ be the subspace of $\ZZ^n_2$ generated by the $\ZZ_2$-characteristic vectors $\{\xi_j : j = 1, \ldots, \ell\}$. So $\ZZ_2(F)$ is a subspace of dimension $\ell$. We will adopt the convention that $\ZZ_2(Q) = 1$. Define an equivalence relation $\sim$ on $E_N$ by
\begin{equation}\label{equ001}
 x \sim y ~ \mbox{if and only if}~ \tau(x) = \tau(y) ~ \mbox{and}~ x=ty ~~ \mbox{for some} ~~ t \in \ZZ_2(F)
\end{equation}
where $F$ is the smallest face containing $\tau(x)$ in its relative interior. The quotient space $N(Q, \xi, \tau) = E_N/\sim$ has a natural $\ZZ_2^n$-action induced from the $\ZZ_2^n$-action on $E_N$. Clearly, the orbit space of $\ZZ_2^n$-action on $N(Q, \xi, \tau)$ is $Q$. Let $$\pi: N(Q, \xi, \tau) \to Q ~ \mbox{defined by} ~ \pi([x]^{\sim}) = \tau(x)$$ be the associated map to the orbit
space $Q$. When $\tau$ is a trivial bundle, we denote $N(Q, \xi, \tau)$ by $N(Q, \xi)$ and $((Q, \xi), (E_N, Q, \tau))$ by $(Q, \xi)$.

The following is proved in Proposition 1.8 of \cite{DJ} when the orbit space $Q$ is a simple convex polytope and in Lemma 3.1 in \cite{LM} when the 2-torus manifold is connected. The similar arguments work to prove the following.
\begin{lemma}\label{clasi}
Let $N$ be a locally standard $2$-torus manifold and $((Q, \xi), (E_N, Q, \tau))$ be the $\ZZ_2$-characteristic model obtained from $N$. Then $N(Q, \xi, \tau)$ is $\ZZ_2^n$-equivariantly homeomorphic to $N$ covering the identity on $Q$.
\end{lemma}

\begin{example}\label{disc1}
Let $Q$ be a closed disc and $\xi : \{\partial{Q}\} \to \ZZ_2^2$ be a map defined by $\xi(\partial{Q}) = (a, b) \in \ZZ_2^2$. Then $\xi$ is a $\ZZ_2$-characteristic function on $Q$. Note that any $\ZZ_2^2$-principal bundle over $Q$ is trivial. So by Lemma \ref{clasi}, $N(Q, \xi)$ is locally standard 2-torus manifold. Since $(a, b) \in \ZZ_2^2$, there is $(c, d) \in \ZZ^2$ such that $ad - bc =1$. So applying an automorphism of $\ZZ_2^2$ we may assume $(c, d) = (1, 0)$ and $(a, b) = (0, 1)$. Then $$N(Q, \xi) \cong (\ZZ_2^2 \times Q)/\sim ~ \cong \ZZ_2 \times (\ZZ_2 \times Q)/\sim ~ \cong \ZZ_2 \times S^2$$ where the action of $\ZZ_2$, determined by $(0, 1)$, on $S^2 \subset \RR^3$ is the reflection along the plane passing through a great circle of $S^2$.
\begin{figure}[ht]
        \centerline{
           \scalebox{0.70}{
            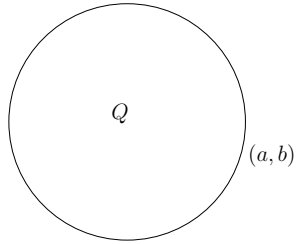
            }
          }
       \caption {A $\ZZ_2$-characteristic function on a closed disc.}
        \label{egsc4}
      \end{figure}
\end{example}

\begin{lemma}
Let $N$ be an $n$-dimensional 2-torus manifold and $((Q, \xi), (E, Q, \tau))$ be the $\ZZ_2$-characteristic data obtained from $N$ such that $\tau$ is trivial. Then $N$ is connected if and only of $\rm{Im}(\xi)$ generates $\ZZ_2^n$.
\end{lemma}

\begin{proof}
Let $F_1, \ldots, F_m$ be the facets of $Q$ and $\xi(F_1), \ldots, \xi(F_m)$ generates $\ZZ_2^n$. Since $E= \ZZ_2^n \times Q$ and $Q$ is connected, $E/\sim$ become connected after applying the identification $\sim$ on $\ZZ_2 \times F_i$ for $i=1, \ldots, m$. On the other hand, let $\xi(F_1), \ldots, \xi(F_m)$ generates a proper submodule $G$ of $\ZZ_2^n$. So $\ZZ_2^n \cong G \oplus \ZZ^k_2$ for some $1 < k < n$. Then $$E/\sim ~~ \cong ((\ZZ_2^k \oplus G)\times Q)/\sim ~~ \cong \ZZ_2^k \times ((G \times Q)/\sim).$$ Therefore this identification space can not be connected.
\end{proof}

\subsection{Orientibility of 2-torus manifolds}\label{ori_2torus}
In this subsection, we give a necessary and sufficient condition for orientability of of 2-torus manifolds generalizing the Theorem 1.7 of Nakayama-Nishimura \cite{NN} and the Theorem 3.1 of Soprunova-Sottile \cite{SS2}. Now we introduce the following definition.

\begin{defn}\label{pol_sp}
An $n$-dimensional nice manifold $Q$ with corners is called a polyhedron with a special vertex if the following holds:
\begin{enumerate}
\item there is a vertex $v$ of $Q$, 

\item there is a subset $Q_v \subset Q$ such that $Q_v$ has a CW-complex structure of dimension $k \leq n-1$ and interior of each $k$ cell belongs to the interior of a $k$-dimensional face of $Q$ not containing the vertex $v$. 

\item $Q \setminus Q_v$ is homeomorphic to $\RR^n_{\geq 0}$ as manifold with corners.
\end{enumerate}
\end{defn}
For simplicity, by a pair $(Q, Q_v)$ we mean a polyhedron with a special vertex where $v$ is a vertex of $Q$ and $Q_v$ is the $k$-dimensional subset of $Q$ which satisfies all the conditions of Definition \ref{pol_sp}.

\begin{example}
\begin{enumerate}
\item All simple polytopes are polyhedron with a special vertex. 

\item Consider the manifold $Q$ with corner of Figure \ref{eqsm02} which is a obtained from 3-cube by deleting the interior of a tetrahedron situated in the interior of the 3-cube. This is a polyhedron with special vertex $v$, since the set $Q_v$ in Figure \ref{eqsm02} satisfy the conditions of Definition \ref{pol_sp}.  

\begin{figure}[ht]
        \centerline{
           \scalebox{.76}{
            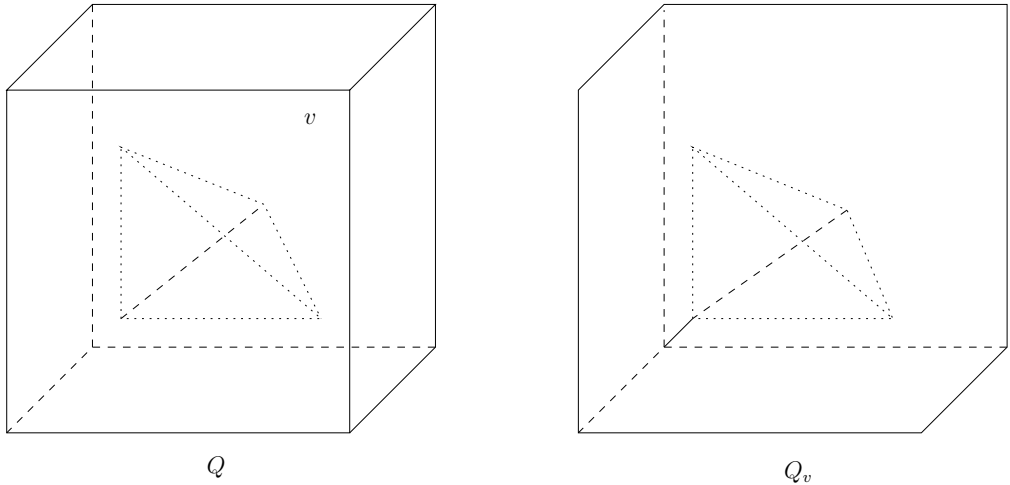
            }
          }
       \caption {A polyhedron with a special vertex.}
        \label{eqsm02}
      \end{figure}
\end{enumerate}
\end{example}

We want to remark that at this moment we do not know if every manifold with corners with a vertex of dimension $n \geq 3$ is a polyhedron with special vertex. Orientability of connected 2-dimensional 2-torus manifolds is remarked at the end of section 4 in \cite{LM}.

\begin{theorem}\label{ori_2tor}
Let $N$ be an $n$-dimensional 2-torus manifold and $((Q, \xi), (E, Q, \tau))$ be the $\ZZ_2$-characteristic data obtained from $N$. Let $\mathcal{F}(Q) = \{F_1, \ldots, F_m\}$, and $(Q, Q_v)$ be a polyhedron with a special vertex where $v \in F_{i_1} \cap \ldots \cap F_{i_n}$, and $$A_i = \{a_{i_j} \neq 0 : \xi_i = a_{i_1} \xi_{i_1} + \cdots + a_{i_n} \xi_{i_n}\}.$$ Then $N$ is orientable if and only if $\dim(Q_v) < n-1$ or the cardinality of $A_i$ is odd for all $i=1, \ldots, m$.
\end{theorem}

\begin{proof}
We prove the theorem using cellular homology. Changing the $\ZZ_2^n$-action on $N$ (if necessary), we may assume that $\xi_{i_j} = e_j$ for $j=1, \ldots, n$ where $e_1, \ldots, e_n$ is the standard basis of $\ZZ_2^n$. Let $Q_v$ be the $k$-dimensional subset of $Q$ which satisfy all the conditions of Definition \ref{pol_sp} and $\{C_1, \ldots, C_s\}$ be the $k$-cells of $Q_v$. Then one can get a CW-complex structure $\{X_i\}_{i=0}^n$ of $N(Q, \xi, \tau) \cong N$ where the $n$-cell is $\pi^{-1}(Q \setminus Q_v) \cong \RR^n$ and other cells come from $\pi^{-1}(Q_v)$. Note that the dimension of cells in $\pi^{-1}(Q_v)$ is $\leq k$. So if $k < n-1$ then $N$ is orientable.

Let $k = n-1$. Then $C_i$ is a subset of the interior of a facet $F_{i_{\ell}}$ for some $i_{\ell} \in \{1, \ldots, m\}$. So $ \pi^{-1}(C_i) = C_i^{\prime} \sqcup C_i^{\prime \prime}$ where the restriction of $\pi$ is a homeomorphism on $C_i^{\prime}$ and $C_i^{\prime \prime}$ to $C_i$. Then the $(n-1)$-cells of the CW-complex structure on $N$ are given by $\sqcup_{i=1}^s (C_i^{\prime} \sqcup C_i^{\prime \prime}) $. Let $Q_C = Q_v \setminus \cup_{i=1}^s C_i$. So $$H_n(X_n, X_{n-1}) \cong H_n(\pi^{-1}(Q)/\pi^{-1}(Q_v)) \cong \ZZ$$ and $$H_{n-1}(X_{n-1}, X_{n-2}) \cong H_{n-1} (\pi^{-1}(Q_v)/\pi^{-1}(Q_C)) \cong \oplus_{i=1}^{2s} \ZZ $$
where each summand is determined by an $(n-1)$-cell of $\sqcup_{i=1}^s (C_i^{\prime} \sqcup C_i^{\prime \prime})$. Note that $C_i^{\prime}$ and $C_i^{\prime \prime}$ are subset of $\partial( \pi^{-1}(Q \setminus Q_v)) \cong S^{n-1}$ and $C_i^{\prime \prime}$ is the homeomorphic image of $C_i^{\prime}$ under $|A_{i_{\ell}}|$ many reflections. Therefore, the degree of the corresponding attaching map is given by $1 + (-1)^{|A_{j_{\ell}}|}$ where $|A_{i_{\ell}}|$ is the cardinality of $A_{i_{\ell}}$. Hence the boundary map $$d_n : H_n(X_n, X_{n-1}) \to H_{n-1}(X_{n-1}, X_{n-2})$$
is zero if and only if $|A_i|$ is odd for $i= 1, \ldots, m$.

\end{proof}

\begin{example}\label{ori_eg}
Consider the $\ZZ_2$-characteristic function given in Figure \ref{eqsm03}. Then $N(Q, \xi)$ is a locally standard 2-torus manifold and it satisfy all the conditions of Theorem \ref{ori_2tor}. So $N(Q, \xi)$ is an orientable $\ZZ_2^3$-manifold.  
\begin{figure}[ht]
        \centerline{
           \scalebox{.76}{
            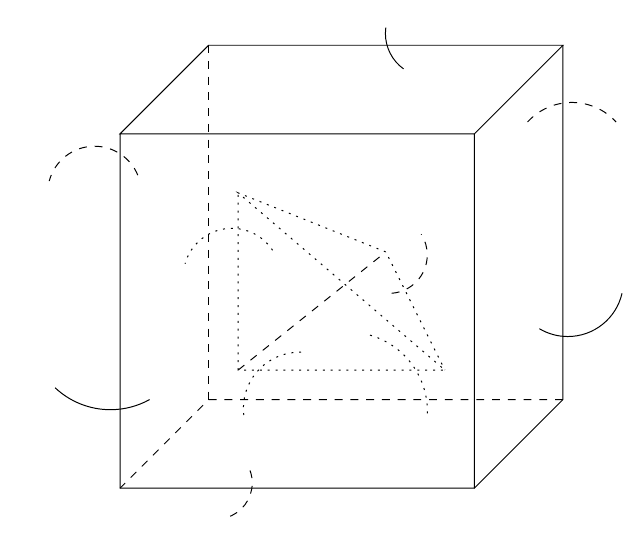
            }
          }
       \caption {A $\ZZ_2$-characteristic model.}
        \label{eqsm03}
      \end{figure}
\end{example}

\begin{remark}
The 2-torus manifold $N(Q, \xi)$ of Example \ref{ori_eg} is not a small cover and not a real torus manifold. So the Theorem \ref{ori_2tor} generalizes the results of Nakayama-Nishimura fo small covers and of Soprunova-Sottile for real toric manifolds. There are some cell complexes, considered on \cite{SS2}, which are not manifold with corners. For example a cone on $m$-gon is a cell complex, but not a manifold with corners if $m > 3$. Note that the manifold with corners of Figure \ref{eqsm03} is not a cell complex as in the Defibition 2.2 of \cite{SS2}. So the orientable spaces determined in \cite{SS2} do not contain our all orientable 2-torus manifolds and vice versa. 
\end{remark}

\subsection{Equivariant connected sum}\label{eqi_con_sum}
This is a well known construction. Here, we briefly discuss equivariant connected sum of $2$-torus manifolds along an orbit. Let $A$ and $B$ be two orbit of $n$-dimensional $2$-torus manifolds $M$ and $N$ respectively. Let the isotropy group of $A$ is isomorphic to the isotropy group of $B$. So $A$ and $B$ are subset of connected component of $M_{i_1} \cap \cdots \cap M_{i_l}$ and $N_{i_1} \cap \cdots \cap N_{i_l}$ for a unique collection of $\ZZ_2$-characteristic submanifolds $\{M_{i_1}, \ldots, M_{i_l}\}$ and $\{N_{i_1}, \cdots, N_{i_l}\}$ of $M$ and $N$ respectively.
Since $\ZZ_2^n$-action is locally standard, there are neighborhood $U_a$ and $U_b$ of $A$ and $B$ such that $U_a$ and $U_b$ are $\ZZ_2^n$-equivariantly homeomorphic to $\cup_{j=1}^{2^{n-{\ell}}} B^{n}$, where $B^n$ is the open $n$-ball in $\RR^n$. Changing the action (if necessary) of $\ZZ_2^n$ on $N$ by an automorphism of $\ZZ_2^n$, we may assume that $\ZZ_2^n$ actions on the $\ZZ_2^n$-invariant neighborhood $U_a$ of $A$ and $U_b$ of $B$ are equivalent. That is we may assume the isotropy of $M_{i_j}$ is same as $N_{i_j}$ for $j=1, \ldots, \ell$. Identifying the boundary  of $M \setminus U_a$ and $ N \setminus U_b$ via an equivariant diffeomorphism we get a smooth manifold, denoted by $M \#_{A, B} N$, with a natural locally standard $\ZZ_2^n$-action. So $M \#_{A, B} N $ is a $2$-torus manifold. Note that if $M$ and $N$ are oriented locally standard $2$-torus manifolds then we can construct oriented connected sum $M \#_{A, B} N$.

Note that an equivariant connected sum of two small covers along a principal orbit is a 2-torus manifold but not a small cover, since the orbit space of the equivariant conncted sum is homeomorphic as manifold with corners to a polytope with holes. Polytope with holes were discussed in section 2 of \cite{PS2}.  

\begin{lemma}\label{eqi_con_tor}
Let $N$ be an oriented 2-torus manifold over a polygon $Q$. Then $N$ is $\ZZ_2^2$-equivariantly homeomorphic to equivariant connected sum of finitely many $T^2$.
\end{lemma}

\begin{proof}
Let $v_0, \ldots, v_m$ be the vertices and $F_i=[v_i, v_{i+1}]$ for $i=1, \ldots, m ~(\mod m+1)$ be the facets (edges) of $Q$. Let $\xi : \mathcal{F}(Q) \to \ZZ_2^2$ be the corresponding $\ZZ_2$-characteristic function and $\xi_i = \xi(F_i)$ for $i=0, \ldots, m$. Since $F_1 \cap F_2 = v_2$ and $\ZZ_2^2$-action on $N$ is locally standard, we may assume $\xi_1 = (1, 0)$ and $\xi_2 = (0, 1)$. Then from Theorem \ref{ori_2tor}, we have $\xi_i=\xi_1 ~\mbox{or}~ \xi_{2}$ for $i=0, \ldots, m$. Since $\{\xi_i, \xi_{i+1}\}$ is a basis for $i = 0, \ldots, m ~ (\mod m+1)$, so $m$ is odd and $\xi_{0} = \xi_2 = \ldots = \xi_{m-1} = (1, 0)$ and $\xi_1 = \xi_3= \ldots = \xi_m = (0, 1)$, see Figure \ref{eqsm01} $(a)$. 

On the other hand, the 2-torus manifold corresponding to the $\ZZ_2$-characteristic model given by Figure \ref{eqsm01} (b) is the torus $T^2$. One can show that the $\ZZ_2$-characteristic function of $\frac{m-1}{2}$ times connected sum of $T^2$ along fixed points is same as $\xi$ on $Q$. For example, a $\ZZ_2$-characteristic function of $T^2 \# T^2$ is given by Figure \ref{eqsm01} (c). Hence by Lemma \ref{clasi}, $N$ is $\ZZ_2^2$-equivariantly homeomorphic to $\frac{m-1}{2}$ times equivariant connected sum of $T^2$. 

\begin{figure}[ht]
        \centerline{
           \scalebox{.72}{
            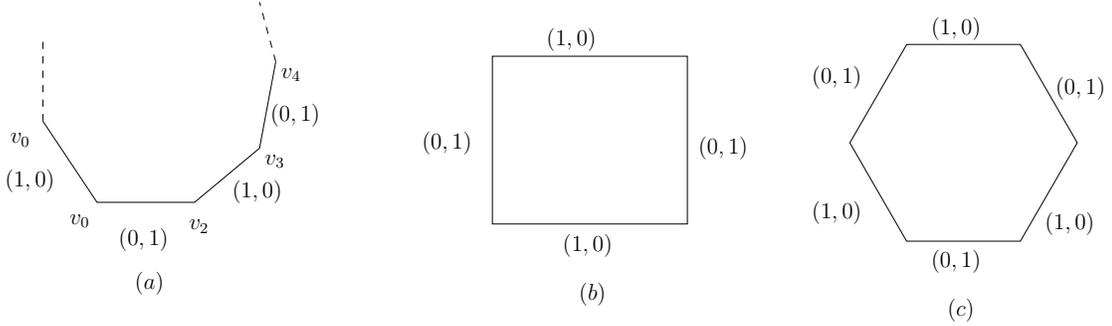
            }
          }
       \caption {Some $\ZZ_2$-characteristic functions.}
        \label{eqsm01}
      \end{figure}

\end{proof}


\section{Construction of manifolds with boundary}\label{def}
In this section we construct $\ZZ_2^n$-manifolds with boundaries where boundaries are locally standard $2$-torus manifolds. 
 
\begin{defn}\label{face-simp}
Let $ Y $ be an $(n+1)$-dimensional compact manifold with corners. $Y$ is said to be {\em face-simple} manifold with exceptional facets if the following holds.
\begin{enumerate}
 \item There exists facets $\{Q_1, \ldots, Q_k\}$ (called exceptional facets) of $Y$ such that $Q_i \cap Q_j$ is empty for $i \neq j$.
\item $V(Y) = \bigcup_{i=1}^k V(Q_i)$.
\item Any codimension-$\ell$ face of $Y$ is the connected component of the intersection of unique $\ell$ many facets.
\end{enumerate}
\end{defn}

Definition \ref{face-simp} generalizes the notion of vertex cut of edge simple polytopes, which was introduced in Section 2 of \cite{Sar2} to study equivariant cobordism quasitoric manifolds. The concept of face-simple polytpe also appears in \cite{Sar4} and \cite{SS}. An example of face-simple manifold with exceptional facets is given in Figure \ref{egc2}. For simplicity, we use the notation $\{Y \backslash Q_1, \ldots, Q_k\}$ to emphasize that $Y$ is a face-simple manifold with exceptional facets $Q_1, \ldots, Q_k$. 
 

Let $\{Y \backslash Q_1, \ldots, Q_k\}$ be an $(n+1)$-dimensional face-simple manifold with exceptional facets and $$ \mathcal{F}(Y) = \{ F_1, \ldots, F_{m}\} \cup \{Q_1, \ldots, Q_k\}.$$ Then we say $F_1, \ldots, F_m$ are the facets of $\{Y \backslash Q_1, \ldots, Q_k\}$ and write as $\mathcal{F}(Y \backslash Q_1, \ldots, Q_k) = \{F_1, \ldots, F_m\}$.
 
\begin{defn}\label{isofun}
A function $ \lambda \colon \{F_1, \ldots, F_{m}\} \to \ZZ_2^{n}$ is called a $\ZZ_2$-isotropy function on $\{Y \backslash Q_1, \ldots, Q_k \}$ if the set of vectors $ \{\lambda(F_{i_1}), \ldots,$ $ \lambda(F_{i_{\ell}})\}$ is a part of a basis of $\ZZ_2^{n}$ whenever the intersection of the facets $ \{F_{i_1}, \ldots,$ $ F_{i_{\ell}}\}$ is nonempty. The vector $\lambda_i :=
\lambda(F_i)$ is called $\ZZ_2$-isotropy vector assigned to the facet $F_i$ for $i=1, \ldots, m$.
\end{defn}
\begin{remark}
The above definition is the generalization of the s-characteristic function of a simple polytope,
see Section 2 of \cite{Sar1}. Let $\{P_{i_1}, \ldots, P_{i_j}\}$ be the facets of $Q_i$ for $i= 1, \ldots, k$. Then there exists a unique facet $F_{i_{\ell}} \in \{F_1, \ldots, F_{m}\}$ such that $P_{i_{\ell}} = Q_i \cap F_{i_{\ell}}$ for $\ell = 1, \ldots j$ and $i = 1, \ldots, k$. We define
$$\lambda^i : \{P_{i_1}, \ldots, P_{i_j}\} \to \ZZ_2^n ~\mbox{by} ~ \lambda^i(P_{i_{\ell}}) = \lambda_{F_{i_{\ell}}} ~\mbox{for} ~ \ell = 1, \ldots, j.$$
Then $\lambda^i$ is a $\ZZ_2$-characteristic function on $Q_i$ for $i=1, \ldots, k$. 
We define $\ZZ_2$-isotropy functions of some face-simple manifold with exceptional facets in Example \ref{egchar2}.
\end{remark}

Let $\eta : E \to Y$ is a principal $\ZZ_2^n$-bundle. These models of $\ZZ_2$-isotropy function and principal bundle will be referred as $\ZZ_2$-isotropy model and denoted by $((Y, \lambda), (E, Y, \eta))$.

Using a $\ZZ_2$-isotropy model $((Y, \lambda), (E, Y, \eta)$, we construct a manifold with boundary. Let $ F $ be a face of $ Y $ of codimension $\ell$, $ 0< \ell \leq n+1$. If $F$ is a face of $Q_i$ for some $i \in \{1, \ldots, k\}$, then $ F $ is the intersection
of a unique collection of $ \ell $ facets $F_{i_1}, F_{i_2}, \ldots, F_{i_{\ell-1}}, Q_i$ of $ Y $. Otherwise, $ F $ is the intersection of a unique collection of $ \ell $ facets $F_{i_1}, F_{i_2}, \ldots, F_{i_{\ell}}$ of $Y$.

Let
\begin{equation}
M(F) = \left\{ \begin{array}{ll} \< \lambda_{i_j} \colon j = 1, \ldots, \ell-1\> \subseteq \ZZ_2^{n} & \mbox{if} ~ F = F_{i_1} \cap \cdots \cap F_{i_{\ell-1}} \cap Q_i\\
\< \lambda_{i_j} \colon j = 1, \ldots, \ell\> \subseteq \ZZ_2^{n} & \mbox{if} ~ F = F_{i_1} \cap \cdots \cap F_{i_{\ell}}
\end{array} \right.
\end{equation}
where $\< \alpha_i : i = 1, \ldots, s \>$ is the subspace  generated by the vectors $\{\alpha_{i} \colon i = 1, \ldots, s\}$.

We will adopt the convention that $M(Y) = 1= M(Q_i)$ for $i=1, \ldots, k$. We define an
equivalence relation $\sim_b$ on $E$ as follows.
\begin{equation}\label{equilam}
x \sim_b y ~ \mbox{if and only if} ~  \eta(x) = \eta(y) ~ \mbox{ and} ~ x = ty ~ \mbox{for some} ~ t \in M(F)
\end{equation}
where $ F $ is the unique face of $ Y $ containing $ \eta(y) $ in its relative interior. We denote the quotient space $ E/ \sim_b $ by $ W(Y, \lambda, \eta)$. The space $W(Y, \lambda, \eta)$ is a $\ZZ_2^n$ space where the action is induced from the  $\ZZ_2^n$-action on $E$. Let $$\mathfrak{q} : W(Y, \lambda, \eta) \to Y$$ be the projection map defined by $\mathfrak{q}([y]^{\sim_b})= \eta(y)$.
\begin{lemma}\label{orbbd}
The space $ W(Y, \lambda, \eta) $ is an $(n+1)$-dimensional $\ZZ_2^n$-manifold with boundary. The boundary is a disjoint union of $2$-torus manifolds.
\end{lemma}

\begin{proof}
 Let $[y]^{\sim_b} \in W(Y, \lambda, \eta)$. So $\eta(y)$ is a unique point in $Y$. We claim that $[y]^{\sim_b}$ has a neighborhood which is $\ZZ_2^n$-equivariantly homeomorphic to $V$ where $V$ is a $\ZZ_2^n$-invariant open subset of $\RR^{n} \times \RR_{\geq 0}$ ($\ZZ_2^n$ action on $\RR^n \times \RR_{\geq 0}$ is induced from the standard action on the first factor). If $\eta(y)$ belongs to the interior of $Y$, there is a neighborhood $U_y$ of $y$ in $Y$ which is homeomorphic
to an $(n+1)$-dimensional open ball in $\RR^{n+1}_{\geq 0}$. Clearly, $$\mathfrak{q}^{-1}(U_y) = \ZZ_2^n \times U_y/ \sim_b ~ = \ZZ_2^n \times U_y \subset \RR^{n} \times \RR_{\geq 0}.$$ So the claim is true in this case.

Let $\eta(y)$ belong to the relative interior of a codimension $\ell (> 0)$ face $F$ of $Y$ where $F$ is not a face of $Q_i$ for any $i \in \{1, \ldots, k\}$. So there is a neighborhood $U_y$ in $Y$ of $\eta(y)$ which is homeomorphic to $\RR^{n+1-\ell} \times \RR^{\ell}_{\geq 0}$ as manifold with corners. The facets of $\RR^{n+1-\ell} \times \RR^{\ell}_{\geq 0}$ are $\{\RR^{n+1-\ell} \times H_i: i =1, \ldots, \ell\}$ where $H_i$ is the facet of $\RR^{\ell}_{\geq 0}$ with $i$-th coordinate zero. Let $F $ be the connected component of $ F_{i_1} \cap \cdots F_{i_{\ell}}$ where $\{F_{i_1}, \ldots, F_{i_{\ell}}\}$ is a  unique collection of $\ell$ many facets of $Y$. By Definition \ref{isofun}, the subgroup $M(F)$ generated by $\{\lambda_{i_1}, \cdots, \lambda_{i_{\ell}}\}$ is a part of a basis of $\ZZ_2^n$.
So $\{\lambda_{i_1}, \cdots, \lambda_{i_{\ell}}\}$ is a basis of $M(F)$ and $\ZZ_2^n \cong \ZZ_2^{n-\ell} \oplus M(F)$.

Suppose the homeomorphism $ \phi \colon U_{y} \to \RR^{n+1-\ell} \times \RR_ {\geq 0} ^{\ell} $ as manifold with corners sends the facet $ F_{i_j}^{\prime} \cap U_{y} $ of $U_y$ to $ \RR^{n+1-\ell} \times H_j $ for all $j = 1, 2, \ldots, \ell$. Define a $\ZZ_2$-isotropy function $ \lambda_{y} $ on the set of all facets of  $ \RR^{n+1-\ell} \times \RR_{\geq 0} ^{\ell} $ by $$ \lambda_{y} (\RR^{n+1-\ell} \times H_j) = \lambda_{i_j} ~~ \mbox{for all} ~~ j = 1, 2, \ldots, \ell.$$ We define an equivalence relation $ \sim_y $ on $ (\ZZ_2^{n} \times \RR^{n+1-\ell} \times \RR_{\geq 0}^{\ell}) $ as follows.
\begin{equation}
(t_1,c_1,x_1) \sim_y (t_2, c_2, x_2) ~ \mbox{if and only if} ~ (c_1,x_1) = (c_2,x_2) ~ \mbox{and}
~ t_1 - t_2 \in M{\phi{(E)}}.
\end{equation}
where $ \phi{(E)} $ is the unique face of $ \RR^{n+1-\ell} \times \RR_{\geq 0} ^{\ell}$ containing $ (c_1, x_1) $ in its relative interior, for a unique face $E$ of $U_y $ and $M{\phi{(E)}} = M(E) $. So for each $ c \in \RR^{n+1-\ell} $ the restriction of $ \lambda_{y} $ on $\{ ( \{ c \} \times H_j : j = 1, 2, \ldots, \ell\} $
defines a $\ZZ_2$-characteristic function (see Definition \ref{z2char}) on the set of facets of
$ \{ c \} \times \RR_{\geq 0}^{\ell}$.
From the constructive definition of small cover given in section 1.5 of \cite{DJ} it is clear that the quotient space $ \{ c \} \times ( M(F) \times \RR_ {\geq 0}^{\ell})/ \sim_y $ is homeomorphic to $ \{c\} \times \RR^{\ell} $. Hence the quotient space
$$ (\ZZ_2^{n} \times \RR^{n+1-\ell} \times \RR_{\geq 0}^{\ell})/\sim_y ~ = ~ \ZZ_2^{n-\ell} \times \RR^{n+1-\ell} \times (M(F) \times \RR_{\geq 0}^{\ell}) / \sim_y ~ \cong ~ \ZZ_2^{n-\ell} \times \RR^{n+1-\ell} \times \RR^{\ell}.$$
Notice that we have the following commutative diagram. 
\begin{equation}
\begin{CD}
(\ZZ_2^{n} \times U_{y}) @>id \times \phi>> (\ZZ_2^{n} \times \RR^{n+1-\ell} \times \RR_{\geq 0}^{\ell}) \\
@V\mathfrak{q} VV  @V \mathfrak{q}_y VV @. \\
(\ZZ_2^{n} \times U_{y}) / \sim @>\phi_y >> (\ZZ_2^{n} \times \RR^{n+1-\ell} \times \RR_{\geq 0}^{\ell})/\sim_y @>\cong>>
\ZZ_2^{n-\ell} \times \RR^{n+1-\ell} \times \RR^{\ell}
\end{CD}
\end{equation}
Since $\mathfrak{q}$ and $\mathfrak{q}_y$ are $\ZZ_2^n$-equivariant quotient maps, and $\phi$ is a homeomorphism as manifold with corners, the commutativity of the diagram  ensures that the lower horizontal map $ \phi_y $ is a $\ZZ_2^n$-equivariant homeomorphism. So the claim is proved for this case.

Now let $\eta(y)$ belongs to the face $Q_j$ for some $j \in \{1, \ldots, k\}$. Since $Y$ is a compact manifold with corners, there exists a color neighborhood $U_{y}$ of $Q_j$ in $Y$. So there is a homeomorphism $g_y : U_{y} \to Q_j \times [0,1)$ as manifold with corners and a $\ZZ_2^n$-equivariant homeomorphism $ \widetilde{g}_y : \eta^{-1}(U_y) \to \eta^{-1}(Q_j) \times [0, 1)$ such that the following diagram commutes, 
\begin{equation}
\begin{CD}
\eta^{-1}(U_y) @>{\widetilde{g}_y}>> \eta^{-1}(Q_j) \times [0, 1)\\
@VV \eta V  @VV \eta \times id V\\
 U_y @>{g_y}>> Q_j \times [0, 1).
\end{CD}
\end{equation} 
Observe that the identification $\sim_b$ in \eqref{equilam} does not affect the second component of $\eta^{-1}(Q_j) \times [0, 1)$ and $\sim_b$ is same as the relation $\sim$ in \eqref{equ001} on $\eta^{-1}(Q_j)$. Let $\{E_{j_1}, \ldots, E_{j_{\ell}}\}$ be the facets of $Q_j$. Then $E_{j_i} = Q_j \cap F_{j_i}$ for a unique facet $F_{j_i}$ of $\{Y \backslash Q_1, \ldots, Q_k\}$. Define $$\xi^j(E_{j_i}) = \lambda(F_{j_i}) \quad \mbox \quad i=1, \ldots, \ell.$$ Then $\xi^j$ is a characteristc function on $Q_j$. Let $E_j = \eta^{-1}(Q_j)$ and $\tau_j : E_j \to Q_j$ be the restriction of $\eta$. So $\tau_j : E_j \to Q_j$ is a principal $\ZZ_2^n$-bundle. So we have the following $\ZZ_2^n$-equivariant homeomorphisms.
\begin{equation}\label{eq002}
(\ZZ_2^n \times U_y)/ \sim_b ~ \cong (\ZZ_2^n \times Q_j)/ \sim_b \times [0, 1) \cong (\ZZ_2^n \times Q_j)/ \sim \times [0, 1)
\cong N(Q_j, \xi^j, \tau_j) \times [0, 1)
\end{equation}
 where $N(Q_j, \xi^j, \tau_j)$ is the locally standard 2-torus manifold associated
to the $\ZZ_2$-characteristic model $((Q_j, \xi^j), (E_j, Q_j, \tau_j))$.


From the above calculation, we get $$W(Y, \lambda, \eta)= \bigcup_{\eta(y) \in Y} \eta^{-1}(U_y)/\sim_b$$ is a compact $\ZZ_2^n$-manifold with boundary where the boundary is the disjoint union of locally standard 2-torus manifolds $\{X(Q_j, \xi^j, \tau_j): j = 1, \ldots, m\}$. 
\end{proof}

We remark that similar construction of manifolds with corners arises in \cite{Sar2, SS3, Sar4, SS}. 
When $\eta$ is trivial we denote the space $W(Y, \lambda, \eta)$ by $W(Y, \lambda)$.

\begin{example}\label{egchar2}
A $\ZZ_2$-isotropy function is given in the Figure \ref{egc2} of the manifold $Y$ with exceptional facets. The exceptional facets are $Q_1, \ldots, Q_5$ where $Q_1,Q_{2}, Q_{3}, Q_{4}$ are triangles
and $Q_{5}$ is a rectangle. The restriction of the isotropy function on $Q_{i}$ gives that the space
$(\ZZ_2^2 \times Q_i)/\sim$ is $\mathbb{RP}^2$ for each $i\in \{1,2,3,4\}$
and $(\ZZ_2^2 \times Q_5)/\sim$ is $T^2$. Hence the space $\sqcup_{1}^4 \mathbb{RP}^2
\sqcup T^2$ is the boundary of $W(Y, \lambda)$.

\begin{figure}[ht]
        \centerline{
           \scalebox{0.60}{
            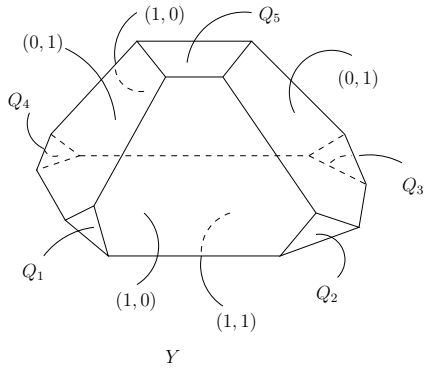
            }
          }
       \caption {An isotropy function of face-simple manifold with exceptional facets.}
        \label{egc2}
      \end{figure}
\end{example}



\section{Equivariant cobordism of 2-torus manifolds}\label{cobor}
Let $\mathfrak{C}$ be the following category: the objects are all oriented $2$-torus manifolds and morphisms are equivariant maps between them. We recall the geometric definition of $\ZZ_2^n$-equivariant cobordism of $2$-torus manifolds. We adhere to the notations of previous Sections. 

\begin{defn}
Two $n$-dimensional oriented $2$-torus manifolds $N_1$ and $N_2$ are said to be $\ZZ_2^n$-cobordant
(or $2$-torus cobordant) if there exists an oriented $\ZZ_2^n$-manifold $W$ with boundary $\partial W $ such that $\partial W $ is $\ZZ_2$-equivariantly homeomorphic to $ N_1 \sqcup - N_2$ under an
orientation preserving homeomorphism. Here $-N_2$ represent the reverse orientation of $N_2$.
\end{defn}
We denote the $\ZZ_2^n$-cobordism class of the oriented $2$-torus manifold $N$ by $[N]_{\delta}$ where $\delta$ represent the $\ZZ_2^n$-action on $N$ or by $[N]$ if the action is clear.
\begin{defn}
The $n$-th $2$-torus cobordism group is the group of all cobordism classes of $n$-dimensional oriented $2$-torus manifolds with the operation of disjoint union. We denote this group by $OS_n$.
\end{defn}


We introduce some sufficient conditions which give affirmative answer of Question \ref{ques}. Let $N$ be an $n$-dimensional 2-torus oriented manifold and $((Q, \xi), (E, Q, \tau))$ be the associated $\ZZ_2$-characteristic model of $N$.  

\begin{theorem}\label{oricob}
Let $\{Y \backslash Q, Q_1, \ldots, Q_k \}$ be an $(n+1)$-dimensional face simple manifold with exceptional facets and $((Y, \lambda), (E_Y, Y, \eta))$ be the associated $\ZZ_2$-isotropy model such that
\begin{enumerate}
\item $\lambda(F_i) = \xi(F_i \cap Q)$ if $F_i \cap Q \neq \emptyset$ for a facet $F_i \in \mathcal{F}(\{Y \backslash Q, Q_1, \ldots, Q_k \})$,

\item $(Y, Y_v)$ is a polyhedron with a special vertex $v \in Q \subset Y$,

\item the bundle $\eta|_{\eta^{-1}(Q)} : \eta^{-1}(Q) \to Q$ is equivariantly isomorphic to $\tau : E \to Q$ and $\tau^i := \eta|_{\eta^{-1}(Q_i)}$ for $i = 1, \ldots, k$.

\item if $v=Q \cap F_{i_1} \cap \cdots \cap F_{i_n}$, $\lambda_i = \lambda(F_i)$ for $F_i \in \mathcal{F}(Y \backslash Q, Q_1, \ldots, Q_k)$ and $A_{F_i}= \{a_{i_j} \neq 0 : \lambda_i = a_{i_1} \lambda_{i_1} + \cdots + a_{i_n} \lambda_{i_n}\}$, then cardinality of $A_{F_i}$ is odd for all $F_i \in \mathcal{F}(Y \backslash Q, Q_1, \ldots, Q_k)$. 
\end{enumerate}
 Then $[N] = [N(Q, \xi, \tau)] = [N(Q_1, \xi^1, \tau^1)] + \cdots + [N(Q_k, \xi^k, \tau_k)]$ in $OS_n$.
\end{theorem}
\begin{proof}
Let $ \mathcal{F}({Y}) := \{ F_{i} \colon i= 1, \ldots, m\} \cup \{ Q_j \colon j = 1, 2, \ldots, k\} $ be the facets of $Y$ and $\lambda_i = \lambda(F_i)$ for $i=1, \ldots, m$. We define a function $ \bar{\xi} \colon \mathcal{F}(Y) \to \ZZ^{n+1}$
as follows.
\begin{equation}\label{bar_char}
 \bar{\xi}( F ) = \left\{ \begin{array}{ll} ( 0, \ldots, 0, 1 ) \in \ZZ^{n+1} & \mbox{if} ~ F = Q_j~ \mbox{and} ~ j~ \in \{1, \ldots, k\} \\
 {(\lambda_i, 0)} \in \ZZ^{n} \times\{0\} \subset \ZZ^{n+1} & \mbox{if}~ F = F_i ~ \mbox{and} ~ i \in \{ 1, \ldots, m\}.
\end{array} \right.
\end{equation}
So the function $ \bar{\xi} $ is a $\ZZ_2$-characteristic function (see Definition \ref{isofun}) on the manifold $Y$ with corners. Let $\bar{\eta}$ be the composition map $$\ZZ_2 \times E_Y \xrightarrow{pr} E_Y \xrightarrow{\eta} Y$$ where $pr$ is the projection map onto a factor. So we get a $\ZZ_2$-characteristic model $((Y, \bar{\xi}), (\ZZ_2 \times E_Y, Y, \bar{\eta}))$. Hence by subsection \ref{condif}, from this $\ZZ_2$-characteristic model we can construct a $2$-torus manifold $ N(Y, \bar{\xi}, \bar{\eta}) $ over $Y$.
There is a natural $ \ZZ_2^{n+1} $ action on $N(Y, \bar{\xi}, \bar{\eta} )$ induced by the $\ZZ_2^{n+1}$-action on $\ZZ_2 \times E_Y$. Let $$\pi : N(Y, \bar{\xi}, \bar{\eta}) \to Y$$ be the orbit map.

Let $e_{n+1} = (0, \ldots, 0, 1) \in \ZZ_2^{n+1}$. We may assume that the vectors $\{e_1, \ldots, e_n\}$ belong to $\ZZ_2^{n+1}$ by making the $(n+1)$-th coordinate zero. From the conditions $(2)$ and $(4)$, and Equation \eqref{bar_char}, one can easily show that the $\ZZ_2$-characteristic model $((Y, \bar{\xi}), (\ZZ_2 \times E_Y, Y, \bar{\eta}))$ satisfies the sufficient conditions of Theorem \ref{ori_2tor}. So $N(Y, \bar{\xi}, \bar{\eta})$ is an $(n+1)$-dimensional oriented 2-torus manifold.

Since $Y$ is a smooth manifold with corners, each $Q_i$ has a collar neighborhood $N_i$ in $Y$ diffeomorphic to $Q_i \ [0, 1)$. Let $A = \cup_{i=1}^k N_i$. Then $\pi^{-1}(Y -A)$ is open subset of $N(Y, \eta)$. So closure of $\pi^{-1}(Y-A)$ is orientable. From Equation \eqref{bar_char} and the construction in Section \ref{def}, we get that the closure of $\pi^{-1}(Y-A)$ in $N(Y, \bar{\xi}, \bar{\eta})$ is nothing but $\ZZ_2 \times W(Y, \lambda, \eta)$ where $W(Y, \lambda, \eta)$
is a $\ZZ_2^n$ manifold with boundary as constructed in Lemma \ref{orbbd}. Hence $W(Y, \lambda, \eta)$ is an oriented manifold with boundary where the boundary is  $N(Q, \xi, \tau) \sqcup N(Q_1, \xi^1, \tau_1) \sqcup \ldots \sqcup N(Q_k, \xi^k, \tau_k)$. This proves the Theorem.
\end{proof}

Equivariant connected sum of locally standard 2-torus manifolds are discussed in Subsection \ref{eqi_con_sum}. Following lemma may be known to experts.

\begin{lemma}\label{coblem}
The equivariant connected sum along the orbits of same isotropy of two oriented locally standard 2-torus manifolds is equivariant cobordant to the disjoint union of these two manifolds.
\end{lemma}
\begin{proof}
Let $N_1$ and $N_2$ be two locally standard 2-torus manifolds of dimension $n$. 
Then $B_1 := [0,1] \times N_1$ and $B_2 := [0,1] \times N_2$ are $\ZZ_2^n$-manifolds with boundary such that $$\partial{B_1}= 0 \times (-N_1) \sqcup 1 \times N_1 ~ \mbox{and} ~ \partial{B_2}= 0 \times (-N_2) \sqcup 1 \times N_2.$$
Let $A_1 \subset N_1$ and $A_2 \subset N_2$ be two orbits of same isotropy group. Let $ U_1 \subset B_1$ and $ U_2 \subset B_2$ be two $\ZZ_2^n$ invariant open neighborhoods of $1\times A_1$ and $1 \times A_2$ respectively. Identifying $\partial{U_1} \cap (B_1 - U_1)$ and $\partial{U_2} \cap (B_2 - U_2)$ via a suitable orientation reversing $\ZZ_2^n$-equivariant homeomorphism we get the lemma.
\end{proof}

We remark that equivariant connected sum of small covers may not be a small cover if the connected sum is performed along a non-trivial orbit.

\begin{corollary}\label{cor_con1}
Let $N$ be an $n$-dimensional oriented locally standard 2-torus manifold and $\ZZ_2^n$-equivariantly homemorphic to equivariant connected sum (not necessarily along fixed points) to small covers $N_1, \ldots, N_k$. Then the cobordism class of $N$ contains a small covers.
\end{corollary}

\begin{proof}
Connected sum of small covers along a fixed point is again a small cover, see \cite{DJ}. So the Lemma \ref{coblem} implies the result.
\end{proof}

\begin{corollary}\label{cor_con2}
If the 2-torus manifolds $N(Q_1, \xi^1, \tau_1), \ldots, N(Q_k, \xi^k, \tau_k)$ in Theorem \ref{oricob} are equivariant connected sum of small covers, then $[N]$ contains a small cover.
\end{corollary}

\begin{proof}
From \cite{DJ}, we know equivariant connected sum of small covers along fixed points again a small cover. Therefore the corollary follows from Theorem \ref{oricob} and Corollary \ref{cor_con1}.
\end{proof}

\begin{prop}
Any oriented 2-dimensional small cover a polygon is $\ZZ^2_2$-equivariantly boundary.
\end{prop}

\begin{proof}
The 2-torus manifold $T^2$ is $\ZZ_2^2$-equivariantly boundary. So the proposition follows from Lemma \ref{eqi_con_tor} and \ref{coblem}.
\end{proof}

\begin{example}
Let $\Delta^{2n} \subset \RR^{2n+1}$ be a $2n$-dimensional simplex with vertices $v_0, \ldots, v_{2n}$ and facets $F_0^{\prime}, \ldots, F_{2n}^{\prime}$ where $F_i^{\prime}$ does not contain the vertex $v_i$ for $i=0, \ldots, 2n$. Cut off a neighborhood of $\{v_0\}, \{v_1\}$ and $F_0 \cap F_1$ by hyperplanes $K_0, K_1$ and $K_2$ respectively such that $K_i \cap K_j \cap \Delta^{2n} = \emptyset$ for $i \neq j$. Then the remaining set $Y$ is a simple $2n$-polytope. Note that $K_i \cap Y = K_i \cap \Delta^{2n}$ is a $(2n-1)$-simplex for $i=0, 1$ and $K_2 \cap Y = K_2 \cap \Delta^{2n} \cong \Delta^1 \times \Delta^{2n-2}$. Let $K_i \cap Y = Q_i$ for $i=0, 1, 2$ and $F_j = F_j^{\prime} \cap Y$ for $j=0, \ldots, 2n$. Then $\{Y \backslash Q_0, Q_1, Q_2\}$ is a face-simple manifold with exceptional facets and  $$\mathcal{F}(\{Y \backslash Q_0, Q_1, Q_2\}) = \{F_0, \ldots, F_{2n}\}.$$ Define a function $\eta : \{F_0, \ldots, F_{2n}\} \to \ZZ_2^{2n-1}$ by $$\eta(F_0) = e_1, \eta(F_{2n})= e_1 + \cdots + e_{2n-1}, ~\mbox{and} ~ \eta(F_i) = e_i ~ \mbox{for} ~ i=1, \ldots, 2n-1.$$ Then $(\{Y \backslash Q_0, Q_1, Q_2\}, \eta)$ is a $\ZZ_2$-isotropy model.

 Let $Q_i : = K_i \cap Y$ for $i=0, 1, 2$. Then the facets of $Q_i$ is given by $$\mathcal{F}(Q_i) = \{F_j \cap Q_i : j \neq i\} ~\mbox{for} ~ i=0, 1 ~ \mbox{and} ~ \mathcal{F}(Q_3) = \{F_j \cap Q_3 : j= 0, \ldots, 2n\}.$$ For each $i \in \{0, 1, 2\}$, define a function $$\xi^i : \mathcal{F}(Q_i) \to \ZZ_2^{2n-1}$$ by $\xi^i(F_j \cap Q_i) = \eta(F_j)$ if $F_j \cap Q_i \neq \emptyset$ and $j \in \{0, \ldots, 2n\}$. Then $(Q_i, \xi^i)$ is a $\ZZ_2$-characteristic model on $Q_i$ for $i=0, 1, 2$. Note that the $\ZZ_2$-isotropy model $(Y, \eta)$ satisfy the conditions of Theorem \ref{oricob}. So $W(Y, \eta)$ is an oriented $\ZZ_2^{2n-1}$-manifold with boundary, where the boundary is $$N(Q_0, \xi^0) \sqcup N(Q_1, \xi^1) \sqcup N(Q_2, \xi^2).$$ Since $Q_i$ is a $(2n-1)$-simplex and the $\ZZ_2$-characteristic function is standard, the space $N(Q_i, \xi^i)$ is $\ZZ_2^{2n-1}$-equivariantly homeomorphic to $\RR P^{2n-1}$ for $i=0, 1$. So by Theorem \ref{oricob}, $N(Q_2, \xi^2)$ is $\ZZ_2^{2n-1}$-equivariantly oriented cobordant to two copies of $\RR P^{2n-1}$. 
\end{example}
%

{\bf Acknowledgement.} 
The author thanks Pacific Institute for Mathematical Sciences, University of Regina and University of Calgary for financial support.

\renewcommand{\refname}{References}

\bibliographystyle{alpha}
\bibliography{bibliography.bib}

\vspace{1cm}

\vfill

\end{document}